\documentclass[10pt]{article}

\usepackage{amssymb,latexsym}
\usepackage{amsfonts}
\usepackage{amsmath,amsthm,graphicx}

\newcounter{cnt1}
\newcounter{cnt2}
\newcounter{cnt3}
\newcommand{\blr}{\begin{list}{$($\roman{cnt1}$)$} {\usecounter{cnt1}
        \setlength{\topsep}{0pt} \setlength{\itemsep}{0pt}}}
\newcommand{\bla}{\begin{list}{$($\alph{cnt2}$)$} {\usecounter{cnt2}
       \setlength{\topsep}{0pt} \setlength{\itemsep}{0pt}}}
\newcommand{\bln}{\begin{list}{$($\arabic{cnt3}$)$} {\usecounter{cnt3}
                \setlength{\topsep}{0pt} \setlength{\itemsep}{0pt}}}
\newcommand{\el}{\end{list}}

\newtheorem{thm}{Theorem}[section]
\newtheorem{lem}[thm]{Lemma}
\newtheorem{prop}[thm]{Proposition}
\newtheorem{defn}[thm]{Definition}
\newtheorem{Exm}[thm]{Example}
\newtheorem{rem}[thm]{Remark}

\newtheorem{cor}[thm]{Corollory}

\newcommand{\ilim}{\mathop{\varprojlim}\limits}

\begin{document}
\title {Versal Deformations of Leibniz Algebras}
\author { Alice Fialowski
\thanks{The work was supported by grants from INSA, the Hungarian Academy of Sciences and by grants OTKA T043641, T043034.}, Ashis Mandal and Goutam Mukherjee$^*$ 
}
\maketitle
\date{}
\noindent
\begin{abstract}
In this work we consider deformations of Leibniz algebras over a field of characteristic zero. The main problem in deformation theory is to describe all non-equivalent deformations of a given object. We give a method to solve this problem completely, namely work out a construction of a versal deformation for a given Leibniz algebra, which induces all non-equivalent deformations and is unique on the infinitesimal level.  
\end{abstract}
\section{Introduction}
Deformations of different algebraic and analytic objects are an important aspect if one studies their properties. They characterize the local behavior in a small neighborhood in the variety of a given type objects.

 After the classical work of Gerstenhaber  in the 60's \cite{G1,G2,G3,G4,G5}, formal deformation theory was generalized in different categories. Computations were made, but the question of obtaining all non-equivalent deformations of a given object was not properly discussed for a long time. The right approach to this is the notion of versal deformation: a deformation which induces all non-equivalent ones. The existence of such a ``versal"  deformation for algebraic categories follows from the work of Schlessinger\cite{Sch}.  For Lie algebras it was worked out in \cite{F}, and one can deduce it to other categories as well. It turns out that (under some minor cohomology restrictions) there exists a versal element, which is universal at the infinitesimal level.
 
  In this work we consider Leibniz algebras and give a construction for a versal element. It is parallel to the general constructions in deformation theory as in \cite{I,P,La,GoM,K}. The first specific method was given for Lie algebras in \cite{FiFu}. Here we are going to work out a similar construction for Leibniz algebras, suitable for explicit computations.
  
The structure of the paper is as follows: In Section $2$ we give the necessary definitions for Leibniz algebras and their cohomology. In Section $3$ we define the notion of deformation over a commutative local algebra base and introduce some basic definitions related to deformations. In Section $4$ we construct a specific  infinitesimal deformation of a Leibniz algebra, which turns out to be universal in the class of infinitesimal deformations. We also introduce the notion of versal deformation: a deformation which is unique on the infinitesimal level and induces all other  deformations. Section $5$ deals with obstructions of extending  a deformation, and in Section $6$ we give the construction of a versal deformation.

\section{Leibniz Algebra and its Cohomology} \label{lcohomology}
Leibniz algebras were introduced by J.L.-Loday \cite{L1,L3} and their cohomology was defined in \cite{LP,L2}. Let us recall some basic definitions.
\begin{defn}
A Leibniz algebra is a $\mathbb{K}$-module $L$, where $\mathbb{K}$ is a field, equipped with a bracket operation that satisfies the Leibniz identity: 
$$[x,[y,z]]= [[x,y],z]-[[x,z],y],~~\mbox{for}~x,~y,~z \in L.$$ 
\end{defn}

Any Lie algebra is automatically a Leibniz algebra, as in the presence of antisymmetry, the Jacobi identity reduces to the Leibniz identity. More  examples of Leibniz algebras were given in \cite {LP}, and recently for instance in \cite{A1, A2}.

Let $L$ be a Leibniz algebra and $M$ a representation of $L$. By definition, $M$ is a $\mathbb{K}$-module equipped with two actions (left and right) of $L$,
$$[-,-]:L\times M\longrightarrow M~~\mbox{and}~[-,-]:M \times L \longrightarrow M ~~\mbox{such that}~$$ 
$$[x,[y,z]]=[[x,y],z]-[[x,z],y]$$
holds, whenever one of the variables is from $M$ and the two others from $L$.

Define $CL^n({L}; {M}):= \mbox{Hom} _\mathbb{K}({L}^{\otimes n}, {M}), ~n\geq 0.$ Let 
$$\delta^n : CL^n({L}; {M})\longrightarrow CL^{n+1}(L; M)$$ 
be a $\mathbb{K}$-homomorphism defined by 
\begin{equation*}
\begin{split}
&\delta^nf(x_1,\cdots,x_{n+1})\\
&:= [x_1,f(x_2,\cdots,x_{n+1})] + \sum_{i=2}^{n+1}(-1)^i[f(x_1,\cdots,\hat{x}_i,\cdots,x_{n+1}),x_i]\\
&+ \sum_{1\leq i<j\leq n+1}(-1)^{j+1}f(x_1,\cdots,x_{i-1},[x_i,x_j],x_{i+1},\cdots,\hat{x}_j,\cdots, x_{n+1}).
\end{split}
\end{equation*}
Then $(CL^*(L; M),\delta)$ is a cochain complex, whose cohomology is called the cohomology of the Leibniz algebra $L$ with coefficients in the representation $M$. The $n$-th cohomology is denoted by $HL^n(L; M)$. In particular, $L$ is a representation of itself with the obvious action given by the bracket in $L$. The  $n$-th cohomology of $L$ with coefficients in itself is denoted by $HL^n(L; L).$

\section{Deformations}
We introduce the notion of deformation of a Leibniz algebra over a commutative algebra base. For an analogous definition for Lie algebras see \cite{F,FiFu}. 
Fix a field $\mathbb{K}$ of characteristic zero. Let $L$ be a Leibniz algebra over $\mathbb{K}$ and $A$ a commutative algebra with identity over $\mathbb{K}$. Let $\varepsilon:{A}\rightarrow
\mathbb{K}$ be a fixed augmentation, that is an algebra homomorphism with $\varepsilon(1)=1$ and $ker (\varepsilon)=\mathfrak{M}$. Throughout the paper we shall assume that $dim(\mathfrak{M}^{k}/\mathfrak{M}^{k+1})<\infty$ for all $k$.
\begin{defn}
A deformation $\lambda$ of ${L}$ with base
$({A},\mathfrak{M})$, or simply with base ${A}$, is an ${A}$-Leibniz
algebra structure on the tensor product
${A}\otimes_\mathbb{K}{L}$ with the bracket $[,]_\lambda$ such that
 \[
 \varepsilon\otimes id:{A}\otimes {L}\rightarrow \mathbb{K}\otimes {L}
 \]
 is a ${A}$-Leibniz algebra  homomorphism (where the $A$-Leibniz algebra structure on $\mathbb{K}\otimes {L}$ is given via $\varepsilon$).
\end{defn}
A deformation of the Leibniz algebra $L$ with base $A$ is called {\it local}
if $A $ is a local algebra over $\mathbb{K}$, and is called
{\it infinitesimal} (or {\it first order}) if, in addition to this, $\mathfrak{M}^2=0$.
Observe that for $l_1,l_2 \in L$ and $a,b \in A$ we have $$[a\otimes l_1,b\otimes l_2]_\lambda = ab[1\otimes l_1,1\otimes l_2]_\lambda$$ by $A$- linearity of $[,]_\lambda$.
Thus to define  a deformation  $\lambda$ it is enough to specify the brackets $[1\otimes l_1,1\otimes l_2]_\lambda$ for $l_1,l_2 \in L$.
Moreover, since $\varepsilon\otimes id:{A}\otimes {L}\rightarrow \mathbb{K}\otimes {L}$
 is a ${A}$-Leibniz algebra  homomorphism, 
 $$(\varepsilon\otimes id)[1\otimes l_1,1\otimes l_2]_\lambda =[l_1,l_2]=(\varepsilon\otimes id)(1\otimes[l_1,l_2])$$
which implies $$ [1\otimes l_1,1\otimes l_2]_\lambda -1\otimes[l_1,l_2] \in ker (\varepsilon \otimes id).$$
Hence we can write 
$$[1\otimes l_1,1\otimes l_2]_\lambda =1\otimes[l_1,l_2]+\sum_{j} c_j \otimes y_j ,$$ where $\sum_{j} c_j \otimes y_j$ is a finite sum with $c_j \in ker(\varepsilon)=\mathfrak{M}$ and $y_j \in L$.
\begin{defn}
Let $A$ be a complete local algebra ( $A=\ilim_{n\rightarrow
\infty}({A}/{\mathfrak{M}^n})$), where $\mathfrak{M}$ is the maximal
ideal in $A$. A formal deformation of $L$ with base $A$ is a $A$-Leibniz algebra structure
on the completed tensor product $$A\hat{\otimes}
L=\ilim_{n\rightarrow
\infty}(({A}/{\mathfrak{M}^n})\otimes L),$$ which is the projective limit of deformations with base ${A}/{\mathfrak{M}^n}$ such that
$$\varepsilon\otimes id:A \hat{\otimes} L \rightarrow K\otimes L=L$$
is a $A$-Leibniz algebra homomorphism.
\end{defn}

\begin{defn}
Suppose $\lambda_1$ and $\lambda_2$ are two deformations of a Leibniz algebra $L$ with base $A$. We call them  equivalent if  there exists a Leibniz algebra isomorphism $$ \phi:(A\otimes L,[,]_{\lambda_1})\rightarrow (A\otimes L,[,]_{\lambda_2})$$ such that $(\varepsilon\otimes id)\circ \phi=\varepsilon\otimes id$.
\end{defn}
We write $\lambda_1 \cong \lambda_2$ if $\lambda_1$ is equivalent to $\lambda_2$. 
\begin{Exm}
If $A= \mathbb{K}[[t]]$ then a formal deformation of a Leibniz algebra $L$ over $A$ is precisely a formal $1$-parameter deformation of $L$(see \cite{B}).
\end{Exm}
\begin{defn}
Suppose $\lambda$ is a given deformation of $L$ with base $(A,\mathfrak{M})$ and augmentation $\varepsilon:{A}\rightarrow \mathbb{K}$. Let $A^\prime$ be another commutative algebra with identity and a fixed augmentation  $\varepsilon^{\prime}:{A^\prime}\rightarrow \mathbb{K}$. Suppose $\phi:A \rightarrow A^{\prime} $ is an algebra homomorphism with $\phi(1)=1$ and $\varepsilon^{\prime} \circ \phi =\varepsilon$. Let $ker(\varepsilon^{\prime})= \mathfrak{M}^\prime$. Then the push-out $\bf{\phi_{*} \lambda}$ is the deformation of $L$ with base $(A^\prime,\mathfrak{M}^\prime)$ and bracket 
          $$[{a_1 }^\prime \otimes_A (a_1\otimes {l_1}),a_2 ^\prime
\otimes_A(a_2\otimes l_2) ]_{\phi_* \lambda}=a_1 ^\prime a_2 ^\prime
\otimes_A[a_1\otimes l_1,a_2\otimes l_2]_\lambda $$
 where $a_1^\prime,a_2 ^\prime \in {A}^\prime,~ a_1,a_2 \in A$ and $l_1,l_2 \in
L$. Here $A^\prime$ is considered as an $A$-module by the map $a^\prime \cdot a=a^\prime \phi(a)$ so that $$A^\prime \otimes L=(A^\prime {\otimes}_{A} A)\otimes L =A^\prime {\otimes}_{A}(A \otimes L).$$
\end{defn}

\begin{rem}\label{push-out}
If the bracket $[,]_\lambda$ is given by 
\begin{equation}\label{exp for bracket}
[1\otimes l_1,1\otimes l_2]_\lambda =1\otimes[l_1,l_2]+\sum_{j} c_j \otimes y_j~\mbox{for}~ c_j \in \mathfrak{M}~\mbox{and}~y_j \in L
\end{equation} 
then the bracket $[,]_{\phi_* \lambda}$ can be written as
\begin{equation}\label{exp for push-out}
[1\otimes l_1,1\otimes l_2]_ {\phi_* \lambda}=1\otimes [l_1,l_2]+\sum_{j}\phi(c_j) \otimes y_j.
\end{equation}
\end{rem}
\section{Universal Infinitesimal Deformation}\label{UID}
In this section we construct a specific infinitesimal deformation of a Leibniz algebra $L$, which turns out to be universal in the class of all infinitesimal deformations of $L$.
Let $L$ be a Leibniz algebra which satisfies the condition $dim (HL^2(L;L)) < \infty$.  This is true for example, if $L$ is finite dimensional. Let us denote the space $HL^2(L;L)$ by $\mathbb{H}$ throughout the paper. Consider the algebra $C_1=\mathbb{K}\oplus \mathbb{H}^\prime$ by setting $(k_1,h_1)\cdot(k_2,h_2)=(k_1 k_2,k_1 h_2+k_2 h_1)$  where $\mathbb{H}^\prime$ is the dual of $\mathbb{H}$ . Observe that the second summand is an ideal of $C_1$ with zero multiplication. Fix a homomorphism $$\mu: \mathbb{H} \longrightarrow
CL^2(L;L)=Hom(L^{\otimes 2};L) $$
which takes a cohomology class into a cocycle representing it. Notice that there is an isomorphism $\mathbb{H}^\prime \otimes L \cong Hom(\mathbb{H}~;L)$, 
so we have
$$C_1 \otimes L 
 =(\mathbb{K} \oplus \mathbb{H}^\prime)\otimes L 
 = (\mathbb{K}\otimes L) \oplus (\mathbb{H}^\prime \otimes L) 
 =L \oplus Hom(\mathbb{H}~;L).$$
Using the above identification, define a Leibniz bracket on $C_1 \otimes L$ as follows.
For $(l_1,\phi_1),(l_2,\phi_2) \in L \oplus Hom(\mathbb{H}~;L)$ let
$$[(l_1,\phi_1),(l_2,\phi_2)]=([l_1,l_2],\psi)$$ where the map  $\psi:\mathbb{H} \longrightarrow L$ is given by 
$$\psi(\alpha)=\mu(\alpha)(l_1,l_2)+[\phi_{1}(\alpha),l_2]+[l_1,\phi_2(\alpha)]~\mbox{for}~\alpha \in \mathbb{H}~.$$ 
It is straightforward to check that $C_1\otimes L$ along with the above bracket is a Leibniz algebra over $C_1$. The Leibniz identity is a consequence of the fact that $\delta \mu(\alpha)=0~ \mbox{for}~ \alpha \in \mathbb{H}$~. Hence we get an infinitesimal deformation of $L$ with base $C_1=\mathbb{K}\oplus \mathbb{H}^\prime $.
\begin{prop}\label{up toisomorphism}
Up to an isomorphism, this deformation does not depend on
the choice of $\mu$. 
\end{prop}
\begin{proof}
Let $$\mu^\prime: \mathbb{H}\longrightarrow CL^2(L;L)$$
be another choice for $\mu$.  
Then for $\alpha\in \mathbb{H}$~, the cochains $\mu(\alpha)$ and
$\mu^\prime(\alpha)$ in $ CL^2(L;L)$
represent the same class $\alpha$. So
$\mu(\alpha)-\mu^\prime(\alpha)$ is a coboundary. Hence we can define a homomorphism 
$$\gamma:\mathbb{H} \longrightarrow CL^1(L;L)$$
on a basis $\{ h_i\}_{1\leq i\leq n}$ of $\mathbb{H}$ by $\gamma(h_i)=\gamma_i$ with $\delta \gamma_i=\mu(h_i)-\mu^\prime(h_i)$. Clearly, $\mu^\prime-\mu=\delta\gamma$. 

Using the identification $C_1\otimes L \cong L \oplus Hom(\mathbb{H}~;L)$, define 
$$\rho:C_1 \otimes L\longrightarrow C_1 \otimes L~~\mbox{by}~
\rho((l,\phi))=(l,\psi),$$ where
$\psi(\alpha)=\phi(\alpha)+\gamma(\alpha)(l)$, $l \in L$ and $\phi \in
Hom(\mathbb{H}~;L)$.

It is routine to check that $\rho$ is a $C_1$-linear automorphism of $C_1\otimes L$, where
 $\rho^{-1}(l,\psi)=(l,\phi)$ with
$\phi(\alpha)=\psi(\alpha)-\gamma(\alpha)(l)$  $\mbox{for}~$ $\alpha \in
\mathbb{H}$~.

It remains to show that $\rho$ preserves the bracket.  Let
$(l_1,\phi_1),(l_2,\phi_2)\in C_1 \otimes L$ with
$\rho(l_1,\phi_1)=(l_1,\psi_1)$ and $\rho(l_2,\phi_2)=(l_2,\psi_2)$.
Suppose 

$[(l_1,\phi_1),(l_2,\phi_2)]=([l_1,l_2],\phi_3)$ 

where $\phi_3(\alpha)=\mu(\alpha)(l_1,l_2)+[\phi_1(\alpha),l_2]+[l_1,\phi_2(\alpha)]$,
and 

$[(l_1,\psi_1),(l_2,\psi_2)]=([l_1,l_2],\psi_3)$

where $\psi_3(\alpha)=\mu^\prime(\alpha)(l_1,l_2)+[\psi_1(\alpha),l_2]+[l_1,\psi_2(\alpha)].$
\begin{equation*}
\begin{split}
\mbox{Then}~\psi_3(\alpha) 
=&\mu^\prime(\alpha)(l_1,l_2)+[\psi_1(\alpha),l_2]+[l_1,\psi_2(\alpha)]\\
=&\mu(\alpha)(l_1,l_2)-\delta\gamma(\alpha)(l_1,l_2)+[\phi_1(\alpha)+\gamma(\alpha)(l_1),l_2]
\\&+[l_1,\phi_2(\alpha)+\gamma(\alpha)(l_2)]\\
=& \mu(\alpha)(l_1,l_2)-[l_1,\gamma(\alpha)(l_2)]-[\gamma(\alpha)(l_1),l_2]
+\gamma(\alpha)([l_1,l_2]) 
\\+&[\phi_1(\alpha),l_2]+[\gamma(\alpha)(l_1),l_2]
+[\phi_1(\alpha),l_2]+[l_1,\phi_2(\alpha)]+[l_1,\gamma(\alpha)(l_2)]\\
=&\mu(\alpha)(l_1,l_2)+[\phi_1(\alpha),l_2]+[l_1,\phi_2(\alpha)]+\gamma(\alpha)([l_1,l_2])\\
=&\phi_3(\alpha)+\gamma(\alpha)([l_1,l_2]).
\end{split}
\end{equation*}
\begin{equation*}
\begin{split}
\mbox{Hence}~\rho([l_1,l_2],\phi_3)=([l_1,l_2],\psi_3)
=[(l_1,\psi_1),(l_2,\psi_2)]
=[\rho(l_1,\phi_1),\rho(l_2,\phi_2)].
\end{split}
\end{equation*}
Therefore, up to an isomorphism, the infinitesimal deformation obtained above is independent of the choice of $\mu$.
\end{proof}
We shall denote this deformation of $L$ by $\eta_1$.
\begin{rem}\label{exp of inf}
Suppose $\{h_i \}_{1\leq i \leq n}$ is a basis of $\mathbb{H}$ and $\{g_i\}_{1\leq i \leq n}$ is the dual basis. Let $\mu(h_i)=\mu_i \in CL^2(L;L)$. Under the identification $C_1 \otimes L = L \oplus Hom(\mathbb{H}~;L)$, an element $(l,\phi)\in L \oplus Hom(\mathbb{H}~;L)$ corresponds to $1\otimes l +\sum_{i=1}^{n}{g_i\otimes \phi(h_i)}$. 
Then for $(l_1,\phi_1),(l_2,\phi_2) \in L \oplus Hom(H;L)$ their bracket $([l_1,l_2],\psi)$ 
corresponds to 
$$1\otimes [l_1,l_2]+ \sum_{i=1}^{n} g_i\otimes (\mu_i(l_1,l_2)+[\phi_1(h_i),l_2]+[l_1,\phi_2(h_i)]).$$
In particular, for $l_1,l_2 \in L$ we have 
$$[1\otimes l_1,1\otimes l_2]_{\eta_1}=1\otimes [l_1,l_2]+\sum_{i=1}^{n}g_i \otimes \mu_i(l_1,l_2).$$ 
\end{rem}
The main property of $\eta_{1}$ is the universality in the class of infinitesimal deformations with a finite dimensional base (Proposition \ref{couniversal}).

Let $\lambda$ be an infinitesimal deformation of the Leibniz algebra $L$ with a finite dimensional base $A$. Let $\{m_i\}_{1\leq i \leq r}$ be a basis of   $\mathfrak{M}=ker(\varepsilon)$ and $\{\xi_{i} \}_{1\leq i \leq r}$ be the dual basis. Note that any element $\xi \in {\mathfrak{M}}^\prime $ can be viewed as an element in the dual space $A^\prime$ with $\xi(1)=0$. For any such $\xi$ set
\begin{equation}\label{alpha}
\alpha_{\lambda,\xi}(l_1,l_2)= \xi \otimes id ([1\otimes l_1,1\otimes l_2]_{\lambda})~~~ \mbox{for}~ l_1,l_2 \in L.
\end{equation}

This defines a cochain $\alpha_{\lambda,\xi} \in Hom(L^{\otimes 2};L)=CL^2(L;L)$.

If we set $\psi_i=\alpha_{\lambda,\xi_i}~\mbox{for}~1\leq i\leq r$, the Leibniz bracket (\ref{exp for bracket}) in terms of the basis of $\mathfrak{M}$ takes the form
\begin{equation}\label{expression}
\begin{split}
[1\otimes l_1,1\otimes l_2]_\lambda 
&=1\otimes[l_1,l_2]+\sum_{i=1}^r m_i \otimes x_i\\
&=1\otimes[l_1,l_2]+\sum_{i=1}^r m_i \otimes \psi_i(l_1,l_2).
\end{split} 
\end{equation}
\begin{lem}
The cochain $\alpha_{\lambda,\xi}  \in CL^2(L;L)$ is a cocycle.
\end{lem}
\begin{proof}
 By definition, 
\begin{equation*}
\begin{split}
&\delta \alpha_{\lambda,\xi}(l_1,l_2,l_3)\\=&[l_1, \alpha_{\lambda,\xi}(l_2,l_3)]+[ \alpha_{\lambda,\xi}(l_1,l_3),l_2]-[ \alpha_{\lambda,\xi}(l_1,l_2),l_3]\\
&- \alpha_{\lambda,\xi}([l_1,l_2],l_3)+ \alpha_{\lambda,\xi}([l_1,l_3],l_2)+ \alpha_{\lambda,\xi}(l_1,[l_2,l_3])~\mbox{for}~l_1,l_2,l_3 \in L.
\end{split}
\end{equation*} 
Observe that 
\begin{equation*}
\begin{split}
&(\xi \otimes id)([1\otimes l_1,[1\otimes l_2,1\otimes l_3]_{\lambda}]_{\lambda})\\
&=(\xi \otimes id)([1\otimes l_1,1\otimes [l_2,l_3]]_{\lambda}+[1\otimes l_1,\sum_{j=1}^{r}m_j\otimes x_j]_{\lambda})~(\mbox{using}~(\ref{expression}))\\
&=\alpha_{\lambda,\xi}(l_1,[l_2,l_3])+\sum_{j=1}^{r}(\xi\otimes id)[1\otimes l_1,m_j \otimes x_j]_{\lambda}.
 \end{split}
\end{equation*}
Moreover,
\begin{equation*}
\begin{split}
(\xi\otimes id)[1\otimes l_1,m_j \otimes x_j]_{\lambda}
&=(\xi\otimes id)m_j[1\otimes l_1,1 \otimes x_j]_{\lambda}\\
&=(\xi\otimes id)m_j(1\otimes [l_1,x_j]+\sum_{i=1}^{r}m_i\otimes x_{ji})\\
&=(\xi\otimes id)(m_j\otimes [l_1,x_j])~~(\mathfrak{M}^2=0)\\
&=[l_1,(\xi\otimes id)(m_j\otimes x_j)].
\end{split}
\end{equation*}
Therefore 
\begin{equation*}
\begin{split}
&(\xi \otimes id)([1\otimes l_1,[1\otimes l_2,1\otimes l_3]_{\lambda}]_{\lambda})\\
&=\alpha_{\lambda,\xi}(l_1,[l_2,l_3])+ [l_1,(\xi\otimes id)\sum_{j=1}^{r}m_j\otimes x_j]\\
&=\alpha_{\lambda,\xi}(l_1,[l_2,l_3])+[l_1,(\xi \otimes id)([1\otimes l_2,1\otimes l_3]_{\lambda}-1\otimes [l_2,l_3])]~~(\mbox{by using}~(\ref{expression}))\\
&=\alpha_{\lambda,\xi}(l_1,[l_2,l_3])+[l_1,\alpha_{\lambda,\xi}(l_2,l_3)]~~( \xi(1)=0).
\end{split}
\end{equation*}
Similarly,$$(\xi\otimes id)([[1\otimes l_1,1\otimes l_2]_\lambda,1\otimes
l_3]_\lambda)
= \alpha_{\lambda,\xi}([l_1,l_2],l_3)+[\alpha_{\lambda,\xi}(l_1,l_2),l_3],$$
$$(\xi\otimes id)([[1\otimes l_1,1\otimes l_3]_\lambda,1\otimes
l_2]_\lambda)=\alpha_{\lambda,\xi}([l_1,l_3],l_2)+[\alpha_{\lambda,\xi}(l_1,l_3),l_2].$$
 It follows that 
\begin{equation*}
\begin{split}
\delta \alpha_{\lambda,\xi}(l_1,l_2,l_3)&= (\xi \otimes id)([1\otimes l_1,[1\otimes
l_2,1\otimes l_3]_\lambda]_\lambda-[[1\otimes l_1,1\otimes
l_2]_\lambda,1\otimes l_3]_\lambda\\
&~~~~~ +[[1\otimes l_1,1\otimes
l_3]_\lambda,1\otimes l_2]_\lambda)\\
&=0 ~~(\mbox{by the Leibniz relation}).
\end{split}
\end{equation*}
\end{proof}
\begin{prop}\label{couniversal}
For any infinitesimal deformation $\lambda$ of a Leibniz algebra
$L$ with a finite dimensional  base $A$ there exists a unique
homomorphism $\phi:C_1=({\mathbb{K}}\oplus
\mathbb{H}^\prime)\longrightarrow A$ such that $\lambda$ is
equivalent to the push-out $\phi_{*}\eta_1$.
\end{prop}
\begin{proof}
Let $\lambda$ be an infinitesimal deformation of a Leibniz
algebra $L$ with base $A$, where $A$ is a finite dimensional local
algebra over $\mathbb{K}$ and $\mathfrak{M}$ is the maximal ideal in $A$.
Let $dim(\mathfrak{M})=r$.
Suppose $\{m_i\}_{1\leq i\leq r}$ is a basis of $\mathfrak{M}$ and $\{\xi_i\}_{1\leq i\leq r}$ be the corresponding dual basis of $\mathfrak{M}^\prime$.
 For $\xi_{i} \in \mathfrak{M}^\prime$ let $a_{\lambda,\xi_{i}} \in
\mathbb{H}$ be the cohomology  class of the cocycle
$\alpha_{\lambda,\xi_{i}}$. The correspondences $$\xi_{i}
\longmapsto \alpha_{\lambda,\xi_{i}}~\mbox{and}~~ \xi_{i} \longmapsto
a_{\lambda,\xi_{i}} $$ for
$1\leq i \leq r$ define homomorphisms
 $$\alpha_\lambda:\mathfrak{M}^\prime \longrightarrow CL^2(L;L)~\mbox{with}~\delta \circ \alpha_\lambda=0~\mbox{and}~ a_\lambda
:\mathfrak{M}^\prime\longrightarrow \mathbb{H}~.$$
We claim that\\
(a)~ Two deformations $\lambda_1$ and $\lambda_2$ are equivalent if
and only if $a_{\lambda_1}=a_{\lambda_2}$.\\
(b)~ If $\phi=id \oplus a^\prime_{\lambda}:C_1 \longrightarrow
\mathbb{K} \oplus \mathfrak{M}=A$ then
$\phi_* \eta_1$ is equivalent to $\lambda$.

Let $\lambda_1$ and $\lambda_2$ be two equivalent deformations of the Leibniz algebra $L$
with base $A$.
Then there exists a $A$-Leibniz algebra isomorphism
$$\rho:(A\otimes L,[,]_{\lambda_1}) \longrightarrow (A\otimes L,[,]_{\lambda_2})
~~\mbox{with}~ (\varepsilon \otimes id)\circ \rho=\varepsilon \otimes id.$$ 
Now
$A \otimes L 
 =(\mathbb{K} \oplus \mathfrak{M})\otimes L
 = (\mathbb{K}\otimes L) \oplus (\mathfrak{M}\otimes L) 
 =L \oplus (\mathfrak{M}\otimes L).$
Thus any element of $A \otimes L$ is of the form $(l,\sum_{i=1}^{r}m_i\otimes l_i)$ where $l_i\in L$ for $1\leq i \leq r$.
By $A$-linearity, $\rho$ is determined by the values $\rho(1\otimes l)~\mbox{for}~l\in L$ and hence $\rho$ is of the form $\rho=\rho_1 +\rho_2$ where $\rho_1:L\longrightarrow L$ and $\rho_2:L\longrightarrow \mathfrak{M}\otimes L$. The map $\rho_1$ must be the identity map  $id:L\longrightarrow L$ by the compatibility $(\varepsilon \otimes id)\circ \rho=\varepsilon \otimes id$.
We shall use the  following standard identifications.
$$Hom(L;\mathfrak{M}\otimes L)\cong \mathfrak{M}\otimes Hom(L;L)\cong Hom(\mathfrak{M}^\prime;Hom(L;L)).$$
In terms of bases of $\mathfrak{M}$ and $\mathfrak{M}^\prime$, the above isomorphisms are given as follows. 
$$\rho_2 \longmapsto \sum_{i=1}^{r}m_i\otimes \phi_i \longmapsto \sum_{i=1}^{r}\chi_i$$ 
where $\phi_i=(\xi_i\otimes id)\circ \rho_2$ and $\chi_i(\xi_j)=\delta_{i,j}\phi_i$. We have
\begin{equation*} 
\begin{split}  
\rho(1\otimes l)
=&\rho_1(1\otimes l)+\rho_2(1\otimes l)
=1\otimes l+\sum_{i=1}^{r}m_i\otimes \phi_i(l)~~\mbox{for}~ l \in L.
\end{split}
\end{equation*}
The map $\rho$ is a Leibniz algebra homomorphism if and only if 
$$\rho([1\otimes l_1,1\otimes l_2]_{\lambda_1})=[\rho(1\otimes l_1),\rho(1\otimes l_2)]_{\lambda_2}.$$
If we take ${\psi_i}^{k}=\alpha_{\lambda_k,\xi_{i}} ~1\leq i \leq r ~\mbox{for}~k=1~\mbox{and}~2$, we have
$$[1\otimes l_1,1\otimes l_2]_{\lambda_k}=1\otimes[l_1,l_2]+ \sum _{i=1}^r {m_i\otimes \psi {_i}^{k}(l_1,l_2)}.$$
\begin{equation*}
\begin{split}
\mbox{Therefore}~&\rho([1\otimes l_1,1\otimes l_2]_{\lambda_1}) \\
=&~ 1\otimes[l_1,l_2]+\sum _{i=1}^{r} m_i \otimes
\phi_i([l_1,l_2])\\
&~+\sum _{i=1}^{r} m_i(1\otimes \psi _{i}^{1}
 (l_1,l_2)
 +\sum_{j=1}^{r} m_j \otimes \phi_j({\psi_{i}}^{1}(l_1,l_2))) \\
=&~ 1\otimes[l_1,l_2]+ \sum _{i=1}^{r} m_i \otimes
\phi_{i}([l_1,l_2])+ \sum _{i=1}^{r} m_i (1 \otimes{
\psi_{i}}^{1} (l_1,l_2))\\
&~~~~~(~ m_i m_j=0).
\end{split}
\end{equation*}
\begin{equation*}
\begin{split}
\mbox{Similarly},~&[\rho( 1\otimes l_1 ),\rho( 1 \otimes l_2 )]_{\lambda_2}\\   
=&~ 1\otimes [l_1,l_2]+\sum _{i=1}^r m_i\otimes\psi_i^2(l_1,l_2)+\sum _{i=1}^r m_i\otimes [l_1,\phi_i(l_2)]\\
&~ + \sum_{i=1}^r m_i\otimes [\phi_i(l_1),l_2].
\end{split}
\end{equation*}
\begin{equation*}
\begin{split}
\mbox{Thus},~&[\rho( 1\otimes l_1 ),\rho( 1 \otimes l_2
)]_{\lambda_2} -\rho([1\otimes l_1,1\otimes l_2]_{\lambda_1}) =0\\
 &\Leftrightarrow \sum _{i=1}^r m_i\otimes
(\psi_i^2(l_1,l_2)-\psi_i^1(l_1,l_2))+\sum _{i=1}^r m_i\otimes
\delta\phi_i(l_1,l_2)=0\\
 & \Leftrightarrow
\psi_i^1(l_1,l_2)-\psi_i^2(l_1,l_2)=\delta\phi_i(l_1,l_2)\\
& ~\mbox{that is},~ \alpha_{\lambda_1,\xi_i}-\alpha_{\lambda_2,\xi_i}=\delta \phi_i ~\mbox{for}~1\leq i \leq r\\
& \Leftrightarrow  a_{\lambda_1} = a_{\lambda_2}.
\end{split}
\end{equation*}
This proves (a).

Now consider the map $$\phi=id \oplus
{a_\lambda}^\prime :C_1 \longrightarrow \mathbb{K} \oplus \mathfrak{M}=A.$$ 
By (a) it is enough to show  that $\alpha_{\phi_{*}\eta_1}=\mu \circ a_\lambda$.
Let $\{ h_i \}_{1\leq i\leq n}$ be a basis of $\mathbb{H}$ and  $\{
g_i\}_{1\leq i\leq n}$ be the corresponding dual basis of $\mathbb{H}^\prime$.  
By Remarks \ref{push-out} and \ref{exp of inf} we have
 $$[1\otimes l_1,1\otimes l_2]_{\phi_{*}\eta_1}=1\otimes [l_1,l_2]+\sum_{i=1}^n
\phi(g_i)\otimes \mu(h_i)(l_1,l_2).$$
Let 
${a_\lambda}^\prime:\mathbb{H}^\prime\longrightarrow \mathfrak{M}$ be the dual of $a_{\lambda}$. Then $$ {a_\lambda}^\prime(g_j)=\sum_{i=1}^r \xi_{i}({a_\lambda}^\prime (g_j))m_i~\mbox{and} ~a_\lambda(\xi_{i})=\sum_{j=1}^n g_j(a_\lambda(\xi_{i}))h_j. $$
\begin{equation*}
\begin{split}
\mbox{Thus}~\alpha_{\phi_{*}\eta_1}(\xi_{i})(l_1,l_2) 
=&~ \xi_{i} \otimes id [1\otimes l_1,1\otimes
l_2]_{\phi_{*}\eta_1}\\
=&~ \xi_{i} \otimes id(1\otimes [l_1,l_2]+\sum_{j=1}^n
\phi(g_j)\otimes \mu(h_j)(l_1,l_2))\\
=&~ \xi_{i} \otimes id(\sum_{j=1}^n {a_\lambda}^\prime(g_j)\otimes
\mu(h_j)(l_1,l_2) )\\
=&~\sum_{j=1}^n \xi_{i} ({a_\lambda}^\prime(g_j))\otimes
\mu(h_j)(l_1,l_2) ) \\
=&~\sum_{j=1}^n g_j(a_\lambda(\xi_{i})) \otimes \mu(h_i)(l_1,l_2) ) \\
=&~\mu(\sum_{j=1}^n g_j(a_\lambda(\xi_{i}))  h_j)(l_1,l_2) \\
=&~\mu \circ a_\lambda(\xi_{i})(l_1,l_2). 
\end{split}
\end{equation*}
The uniqueness part of the theorem follows from the definition
of $\phi$. 
\end{proof}

Suppose $A$ is a local algebra with the unique maximal ideal $\mathfrak{M}$ and $\pi:A\rightarrow A/{\mathfrak{M}^{2}}$ the corresponding quotient map. Assume  $dim(A/{\mathfrak{M}^{2}})$ is finite. The algebra $A/{\mathfrak{M}^{2}}$ is obviously local with maximal ideal ${\mathfrak{M}}/{\mathfrak{M}^{2}}$ and $({\mathfrak{M}}/{\mathfrak{M}^{2}})^{2}=0$. If $\lambda$ is a deformation of $L$ with base $A$ then $\pi_{*} \lambda$ is a deformation with base $A/{\mathfrak{M}^{2}}$ and it is clearly infinitesimal. Therefore by the previous proposition we have a map 
$$a_{\pi *\lambda}:({\mathfrak{M}}/{\mathfrak{M}^{2}})^\prime \rightarrow \mathbb{H}~. $$
\begin{defn}

The dual space $({\mathfrak{M}}/{\mathfrak{M}^{2}})^\prime$  is called the tangent space of A and is denoted by $TA$. The map $a_{\pi *\lambda}$ is called the differential of $\lambda$ and is denoted by $d{\lambda}$.
\end{defn}
It follows from Proposition \ref{couniversal} that equivalent deformations have the same differential.
We have constructed in this section the universal infinitesimal deformation and our goal is to extend it to a versal one. It is known that in the category of deformations of an algebraic object generally there is no universal deformation \cite{Hart}. But under certain natural conditions it is possible to get a ``versal'' object, which still induces all non-equivalent deformations.  

\begin{defn}
A formal deformation $\eta$ of a Leibniz algebra $L$ with base $C$ is called versal, if\\
(i)~for any formal deformation $\lambda$ of $L$ with  base $A$ there exists a homomorphism $f:C \rightarrow A$ such that the deformation $\lambda$ is equivalent to $f_{*}\eta$; \\
(ii)~if $A$ satisfies the condition ${\mathfrak{M}}^2=0$, then $f$ is unique. 
\end{defn}
\begin{thm}
If $\mathbb{H}$ is finite dimensional, then there exists a versal deformation.
\end{thm}
\begin{proof}
Follows from the general theorem of Schlessinger \cite{Sch}, like it was shown for Lie algebras in \cite{F}.
\end{proof}

\section{Obstructions}\label{OBS}
The aim of this section is to study obstructions in extending deformations. For this we need the interpretation of $1$- and $2$-dimensional Harrison cohomology of a commutative algebra. Let us recall some definitions and results from \cite{H}.

Let $A$ be a commutative algebra with $1$ over $\mathbb{K}$ . Let $(C_q(A), \delta)$ denote the standard Hochschild complex, where $C_q(A)$ is the $A$-module $A^{ \otimes (q+1)}$ with $A$ acting on the first factor by multiplication of $A$. Let $Sh_q(A)$ be the $A$-submodule of $C_q(A)$ generated by chains 
\begin{equation*}
\begin{split}
&s_p(a_1,a_2,\ldots,a_q)\\
=&~\sum_{(i_1,i_2,\ldots,i_q)\in Sh(p,q-p)}    {sgn(i_1,i_2,\ldots,i_q)(a_{i_1},a_{i_2},\ldots,a_{i_q}) \in C_q (A)}
\end{split}
\end{equation*}
$\mbox{for}~ a_1,a_2,\ldots,a_q \in A  ~;~0<p<q $.

It turns out that $Sh_{*}$ is a subcomplex of $C_{*}(A)$ and hence we have a complex called the {\it Harrison complex} 
\[Ch(A)=\{ Ch_q(A),\delta\}~;~ Ch_q(A)=C_q(A)/Sh_q(A).\]
For an $A$-module $M$, the Harrison cochain complex defining the Harrison cohomology with coefficients in $M$ is given by
$Ch^{q}(A~;M)=Hom(Ch_{q}(A), M)$.
\begin{defn}
 For an $A$-module $M$ we define
$$ H_{Harr}^{q}(A;~M)= H^{q}(Hom (Ch(A),M)).$$
\end{defn}

\begin{prop}\label{coefficients in M}
Let $A$ be a commutative local algebra with the maximal ideal $\mathfrak{M}$, and let $M$ be an $A$-module with $\mathfrak{M}M=0$.  Then we have the canonical isomorphism
$$ H_{Harr}^{q}(A;~M)\cong  H_{Harr}^{q}(A;~\mathbb{K})\otimes M.$$ 
\end{prop}
\begin{defn}
An extension $B$ of an algebra $A$ by an $A$-module $M$ is a $\mathbb{K}$-algebra $B$ together with an exact  sequence of $\mathbb{K}$-modules $$0\longrightarrow
{M}\stackrel{i}{\longrightarrow} B
\stackrel{p}{\longrightarrow}A\longrightarrow 0,$$ where $p$ is an $\mathbb{K}$-algebra homomorphism, and the $B$-module structure on $i(M)$ is given by the $A$-module structure of $M$ by $i(m) \cdot b=i(m p(b))$.
\end{defn}

\begin{prop}\label{cohomology class corresponds to extension}
\begin{description}
 \item[({\it i})] The space $ H_{Harr} ^1 (A ;M)$ is isomorphic to the space of derivations $A \longrightarrow M$.
\item[({\it ii})] Elements of $ H_{Harr} ^2 (A;M)$  correspond bijectively to isomorphism classes of extensions 
$$ 0\longrightarrow M\longrightarrow B \longrightarrow A\longrightarrow 0 $$
of the algebra $A$ by means of $M$.
\item[({\it iii})]The space $H_{Harr}^1(A;M)$ can also be interpreted as the group of automorphisms of any given  extension
 of $A$ by $M$.
\end{description}
\end{prop}

\begin{cor}\label{TA}
 If $A$ is a local algebra with the maximal ideal $\mathfrak{M}$,  then $$H_{Harr} ^1(A~;\mathbb{K})\cong(\frac{\mathfrak{M}}{\mathfrak{M}^2})'=TA.$$
\end{cor}
Let $\lambda$ be a deformation of a Leibniz algebra $L$ with a finite dimensional local base $A$ and augmentation $\varepsilon$. Consider $[f] \in H^2_{Harr}(A~;\mathbb{K})$. Suppose 
$$0\longrightarrow
\mathbb{K}\stackrel{i}{\longrightarrow} B
\stackrel{p}{\longrightarrow}A\longrightarrow 0$$ is a representative of  the class of $1$- dimensional extensions of $A$, corresponding to the cohomology class of $f$.  Let 
$
I=(i\otimes id):L \cong \mathbb{K}\otimes L \longrightarrow B\otimes L ,\\
P=(p\otimes id): B\otimes L \longrightarrow A\otimes L ~\mbox{and}~
E=(\hat{\varepsilon}\otimes id): B\otimes L \longrightarrow \mathbb{K}\otimes L \cong L,	
$
where $\hat{\varepsilon}=\varepsilon \circ p$ is the augmentation of $B$ corresponding to the augmentation $\varepsilon$ of $A$. Fix a section $q:A \longrightarrow B$ of $p$ in the above extension, then 
\begin{equation}\label{induced by the section} 
b\longmapsto (p(b),i^{-1}(b-q\circ p(b)))
\end{equation} 
is a $\mathbb{K}$ - module isomorphism  $B \longrightarrow (A \oplus \mathbb{K})$. 
Let us denote by $(a,k)_q \in B$ the inverse of $(a,k)\in (A \oplus \mathbb{K})$ under the above isomorphism. The algebra structure of $B$ is determined by $f$ and is given by 
\begin{equation}\label{algebra structure}
(a_1,k_1)_q \circ (a_2,k_2)_q =(a_1 a_2~ ,~ a_1 \cdot k_2+a_2 \cdot k_1+f(a_1,a_2))_q.
\end{equation}
Suppose $dim(A)=r+1$ and $\{m_i\}_{1\leq i \leq r}$ is a basis of the maximal ideal $\mathfrak{M}_A$ of A.  Then $\{n_i \}_{1\leq i \leq {r+1}}$ is a basis of the maximal ideal $\mathfrak{M}_B =p^{-1}(\mathfrak{M}_A)$ of $B$, where $n_j=(m_j,0)_q,~\mbox{for} ~1\leq j\leq r$ and $n_{r+1}=(0,1)_q$. Take the dual basis  $\{\xi_{i}\}_{1\leq i \leq r}$ of ${\mathfrak{M}^\prime _A}$ . Then by (\ref{alpha}) and (\ref{expression}), we have 2-cochains $\psi_i=\alpha_{\lambda,\xi_i} \in CL^2(L;L)$ for $1\leq i \leq r$ such that $[,]_{\lambda}$ can be written as 
$$
[1\otimes l_1,1\otimes l_2]_{\lambda}=1\otimes [l_1,l_2]+\sum_{i=1}^{r} m_i \otimes \psi_i(l_1,l_2)~\mbox{for}~ l_1,l_2 \in L.$$ 
Let $\psi \in CL^2(L;L)=Hom(L^{\otimes 2};L)$ be an arbitrary element. Define a $B$-bilinear operation 
$(B\otimes L)^{\otimes 2} \longrightarrow B\otimes L$,
$$\{b_1 \otimes l_1,b_2 \otimes l_2\}=b_1b_2 \otimes [l_1,l_2]+\sum_{j=1}^{r}b_1b_2n_j \otimes \psi_j(l_1,l_2) +b_1b_2n_{r+1}\psi(l_1,l_2).$$
It is routine to check that the $B$-bilinear map $\{,\}$ satisfies 
\begin{equation}\label{condition for extension}
\begin{split}
&(i)P\{l_1,l_2 \}=[P(l_1),P(l_2)]_\lambda ~~\mbox{for}~l_1,l_2 \in
B\otimes L\\
&(ii)\{I(l),l_1\}=I[l,E(l_1)]~~\mbox{for}~l \in L~\mbox{and}~l_1 \in
B\otimes L.
\end{split}
\end{equation}
So the Leibniz algebra structure $\lambda$ on $A \otimes L$ can be lifted to a $B$-bilinear operation  
$\{ ,\}:(B\otimes L)^{\otimes 2} \longrightarrow B\otimes L$ satisfying (\ref{condition for extension}).

Define $\phi:(B\otimes L)^{\otimes 3}\longrightarrow B\otimes L~~\mbox{by}$
\begin{equation}\label{phi} 
 \phi(l_1,l_2,l_3)=\{l_1,\{l_2,l_3\}\}-\{\{l_1,l_2\},l_3\}+\{\{l_1,l_3\},l_2
\} ~~\mbox{for}~ l_1,l_2,l_3 \in B \otimes L.
\end{equation}
It is clear that $\{,\}$ satisfies the Leibniz relation if and only if $\phi =0 $.
Now from property ({\it i}) in (\ref{condition for extension}) and the definition of $\phi$ it follows that $P\circ \phi(l_1,l_2,l_3)=0~~\mbox{for}~l_1,l_2,l_3 \in B\otimes L$.
Therefore $\phi$ takes values in $ker(P)$.
Observe that $\phi(l_1,l_2,l_3)=0$, whenever one of the arguments belongs to $ker(E)$. 
Suppose $l_1=b\otimes l \in ker (E)\subseteq B\otimes L$.
Since $ker(E)=ker(\hat{\varepsilon})\otimes L=p^{-1}(ker(\varepsilon))\otimes L={\mathfrak{M}}_B \otimes L$, we can write $l_1=\sum_{j=1}^{r+1}n_j \otimes {l_j}^\prime$ with ${l_j}^\prime \in L;~ j=1,\ldots,r+1$. Then for $l_2,l_3 \in B \otimes L$, we get
$$\phi(l_1,l_2,l_3)=\phi(\sum_{j=1}^{r+1}n_j \otimes l^\prime_j,l_2,l_3)
=\sum_{j=1}^{r+1}n_j \phi(l^\prime_j,l_2,l_3)=0.$$
This is because $\phi(l^\prime_j,l_2,l_3) \in ker(P)=im(I)=im(i)\otimes L=i(\mathbb K)\otimes L$ and for any element $k \in \mathbb{K}$ and $l \in L$, 
\begin{equation*}
\begin{split}
&n_j \cdot i(k)\otimes l= i(p(n_j) k)\otimes l= i(m_j \cdot k)\otimes l=i(\varepsilon(m_j) k)\otimes l=0~~\mbox{for}~1\leq j\leq r\\
&\mbox{and}~n_{r+1} \cdot i(k)\otimes l =k n^2_{r+1}\otimes l=0 ~~(m_j \in \mathfrak{M}\subset A~\mbox{and}~ m_j\cdot k=\varepsilon(m_j)k).
\end{split}
\end{equation*} 
The other two cases are similar. Thus $\phi$  defines a linear map $$ \tilde{\phi}:(\frac{B\otimes
L}{ker(E)})^{\otimes 3}\longrightarrow ~~ker(P), $$
$\tilde{\phi}(b_1 \otimes l_1+ker(E),b_2 \otimes l_2+ker(E),b_3 \otimes l_3+ker(E))=\phi(b_1 \otimes l_1,b_2 \otimes l_2,b_3 \otimes l_3)$. 
Moreover, the surjective map $E: B \otimes L \longrightarrow \mathbb{K}\otimes L {\cong}L$,   defined by $b \otimes l \longmapsto \hat{\varepsilon}(b) \otimes l$, induces an isomorphism $\frac{B\otimes
L}{ker(E)} \stackrel{\alpha}{\cong} L$,  where $$\alpha :L \longrightarrow \frac{B\otimes
L}{ker(E)}  ~~;~~\alpha(l)=1 \otimes l + ker(E).$$
Also, $ker(P)=im(I)=i(\mathbb{K})\otimes L =\mathbb{K}~i(1)\otimes L \stackrel{\beta}{\cong} L$ where the isomorphism $\beta$ is given by $\beta(k n_{r+1}\otimes l) =k l$ with inverse ${\beta}^{-1}(l)=n_{r+1}\otimes l$.
Thus we get a linear map $\bar{\phi}:L^{\otimes 3}\longrightarrow L$, such that 
$\bar{\phi}=\beta \circ \tilde{\phi}\circ {\alpha}^{\otimes 3}.$ 
The cochains $\bar{\phi} \in CL^3(L;L)$ and $\phi$ are related by
\begin{equation}\label{relation of phi}
n_{r+1}\otimes \bar{\phi}(l_1,l_2,l_3)= \phi(1\otimes l_1,1\otimes l_2,1\otimes l_3).
\end{equation}
We claim that the cochain $\bar{\phi}$ is a cocycle.
The coboundary $\delta \bar{\phi}$ consists of $10$ terms . Let us rewrite the first term of  $\beta^{-1}\circ \delta \bar{\phi}$ as follows.
\begin{equation*}
\begin{split}
&\beta^{-1}([l_1,\bar{\phi}(l_2,l_3,l_4)])\\
=&~ n_{r+1}\otimes [l_1,\bar{\phi}(l_2,l_3,l_4)]\\
=&~ I([l_1,\bar{\phi}(l_2,l_3,l_4)])~~~(i(1)=n_{r+1})\\
=&~I([l_1,E(1 \otimes \bar{\phi}(l_2,l_3,l_4))])\\
=&~ \{I(l_1),1 \otimes \bar{\phi}(l_2,l_3,l_4)\}~~~(\mbox{by}~(\it ii)~\mbox{of}~ (\ref{condition for extension}))\\
=&~ \{n_{r+1}\otimes l_1,1 \otimes \bar{\phi}(l_2,l_3,l_4)\}\\
=&~ \{1 \otimes l_1,n_{r+1} \otimes \bar{\phi}(l_2,l_3,l_4)\}\\
=&~ \{1 \otimes l_1,\phi(1 \otimes l_2,1 \otimes l_3,1 \otimes l_4)\}~~~(\mbox{by}~(\ref{relation of phi}))\\
=&~ \{1 \otimes l_1,\{1 \otimes l_2,\{1 \otimes l_3,1 \otimes l_4\}\} \} -\{1 \otimes l_1,\{\{1 \otimes l_2,1 \otimes l_3\},1 \otimes l_4\} \} \\
&~ + \{1 \otimes l_1,\{\{1 \otimes l_2,1 \otimes l_4\},1 \otimes l_3\}\}.
\end{split}
\end{equation*}
Similarly, computing other terms and substituting in the expression of $\beta^{-1}\circ \delta\bar{\phi}$, we get
 $\beta^{-1}\circ \delta\bar{\phi}(l_1,l_2,l_3,l_4)=0.$
 
Let us show now that the cohomology class of $\bar{\phi}$ is independent of the choice of the lifting $\{,\}$. Suppose $\{ ,\}$ and $\{ ,\}^\prime$ are two $B$-bilinear operations on $B \otimes L$, lifting the Leibniz algebra structure $\lambda$ on $A \otimes L$.
Let $\bar{\phi}~~\mbox{and}~\bar{\phi'}$ be the corresponding cocycles. Set $\rho =\{,\}^\prime -\{,\} $. 
Then $\rho:(B \otimes L)^{\otimes 2} \longrightarrow B \otimes L
~~\mbox{is a}~B\mbox{-linear map}.$ Observe that
\begin{equation*}
\begin{split}
P\circ \rho(l_1,l_2)
= [P(l_1),P(l_2)]_{\lambda}- [P(l_1),P(l_2)]_{\lambda}=0~~
(\mbox{by}~(\it i)~\mbox{in}~(\ref{condition for extension}) ).
\end{split}
\end{equation*}
Thus $\rho$ takes values in $ker(P)$ and 
induces a linear map 
$$\tilde {\rho}:({\frac{B\otimes
L}{ker(E)}})^{\otimes 2} \longrightarrow ker(P),$$
$$\tilde {\rho}(l_1+ker(E),l_2+ker(E))={\rho}(l_1,l_2)~\mbox{for}~l_1,l_2 \in B \otimes L.$$
Hence we get a $2$-cochain $\bar{\rho}:L^{\otimes 2}\longrightarrow L$ such that 
$\bar{\rho}=\beta \circ \tilde{\rho}\circ \alpha^{\otimes 2} \in CL^2(L;L)$. As before, for $l_1,l_2 \in L$, we have
$n_{r+1} \otimes \bar{\rho}(l_1,l_2) ={\rho}(1 \otimes l_1,1 \otimes l_2).$ Then a straightforward computation yields
 $$\beta^{-1} \circ \delta \bar{\rho}(l_1,l_2,l_3)=\beta^{-1} \circ (\bar{\phi}^\prime -\bar{\phi})(l_1,l_2,l_3),~~\mbox{for}~l_1,l_2,l_3 \in L .$$
Hence $(\bar{\phi}^\prime -\bar{\phi})=\delta \bar{\rho}$.

Suppose a $B$-bilinear operation $\{,\}$ is given on $B \otimes L$, lifting the Leibniz algebra structure $[,]_{\lambda}$ on $A\otimes L$. Then any other $B$-bilinear operation on $B \otimes L$, lifting $[,]_{\lambda}$, is determined by a $2$-cochain $\rho$ as follows. Define  
$\{ ,\}^\prime:(B\otimes L)^{\otimes 2} \longrightarrow B\otimes L$ by
$ \{l_1 ,l_2\}^\prime=\{l_1 ,l_2\}+I\circ \rho(E(l_1),E(l_2))$ $ \mbox{for}~ l_1,l_2 \in B \otimes L$. Then it is easy to see that $\{,\}^\prime$ is a lifting of $[,]_{\lambda}$ such that the $2$-cochain $\bar{\rho}$ induced by the difference $\{ ,\}^\prime -\{ ,\}$, is the given $2$-cochain $\rho$.

The above consideration defines a map 
$\theta_{\lambda}:  H_{Harr}^2(A;\mathbb{K}) \longrightarrow  HL^3(L;L)$ $\mbox{by}~
\theta_{\lambda}([f])=[\bar{\phi}]$, where $[\bar{\phi}]$ is the cohomology class of $\bar{\phi}$. The map $\theta_{\lambda}$ is called the {\it obstruction map}.
\begin{prop}
Let $\lambda$ be a deformation of the Leibniz algebra $L$ with base $A$ and let $B$ be a $1$-dimensional extension of $A$ corresponding to the cohomology class $[f] \in H^{2}_{Harr}(A;\mathbb{K})$. Then $\lambda$ can be extended to a deformation of $L$ with base $B$ if and only if the obstruction $\theta_{\lambda}([f])=0$.
\end {prop}
\begin{proof}
Suppose  $\theta_{\lambda}([f])=0$.  Let 
\begin{equation}\label{1-dimensional ext}
0\longrightarrow
\mathbb{K}\stackrel{i}{\longrightarrow} B
\stackrel{p}{\longrightarrow}A\longrightarrow 0
\end{equation}
 be a $1$-dimensional extension representing the cohomology class $[f]$. 
Let $\{,\}$ be a lifting of the Leibniz algebra structure $\lambda$ on $A \otimes L$ to a $B$-bilinear operation on $B\otimes L$. Let  
$\bar{\phi}$ be the associated cocycle  in $CL^3(L;L)$ as described above.
Then $\theta_{\lambda}([f])=[\bar{\phi}]=0$ implies 
$\bar{\phi}=\delta \rho$ for some $\rho \in CL^2(L;L)$.
Now take $\rho^\prime =-\rho$, and  define a new linear map
$$\{ ,\}^\prime:(B\otimes L)^{\otimes 2} \longrightarrow B\otimes L ~~\mbox{by}~ \{l_1 ,l_2\}^\prime=\{l_1 ,l_2\}+I\circ \rho^\prime(E(l_1),E(l_2)). $$
If $\bar{\phi'}$ denotes the cocycle corresponding to $\{,\}'$, we have $\bar{\phi'}-\bar{\phi}= \delta \bar{\rho'}=-\bar{\phi}$. Hence $\bar{\phi'}=0$ which implies $\phi'=0$. Therefore, $\{,\}'$ is a Leibniz algebra structure on $B \otimes L$ extending $\lambda$. The converse is clear.
\end{proof}
Assume $\theta_{\lambda}([f])=0$ for $[f] \in H^{2}_{Harr}(A;\mathbb{K})$.
Let us denote by $S$ the set of all isomorphism classes of deformations $\mu$ of $L$ with base $B$ such that $p_{*}\mu=\lambda$.
The group of automorphisms $\mathcal{A}$ of the extension (\ref{1-dimensional ext})
has a natural action $\sigma_1$ of $\mathcal{A}$ on $S$, given by $\mu \mapsto u_{*}\mu$ for $u \in \mathcal{A}$. This is clearly well-defined, because if $\mu \cong \mu'$, then $u_{*}\mu \cong u_{*}\mu^\prime$.

On the other hand, $\mathbb{H}$ acts on $S$ as follows. Suppose $\mu$ and $\mu'$ are two deformations of $L$ with base $B$ such that $p_{*}\mu=p_{*}\mu'=\lambda$. Let $\bar{\phi}_\mu$ be the $3$-cochain determined by $[,]_{\mu}$ as before. If $\psi \in CL^2(L;L)$ is the cochain determined by $[,]_\mu -[,]_{\mu^\prime}$, then we have, $\bar{\phi}_{\mu}-\bar{\phi}_{\mu^\prime}=\delta \psi$. But $\bar{\phi}_{\mu}=0=\bar{\phi}_{\mu^\prime}$  ($\mu, \mu'$ are the  Leibniz brackets). Hence $\psi$ is a cocycle.

Suppose now $\mu_1 \cong \mu$ with isomorphism $\rho:(B\otimes L, [,]_{\mu})\longrightarrow (B\otimes L, [,]_{\mu_1})$. As explained in the proof of Proposition \ref{couniversal}, the $2$-cochain determined by $[,]_{\mu}-[,]_{\mu_1}$ is a coboundary. Hence $[,]_\mu -[,]_{\mu_1}=\delta b_\rho$ for some $1$-cochain $b_\rho$, determined by $\rho$. Therefore, 
$ \psi_1=[,]_{\mu'}-[,]_{\mu_1}=([,]_{\mu'}-[,]_{\mu})+([,]_{\mu}-[,]_{\mu_1})=\psi +\delta b_{\rho}.$
This allows us to introduce a map $\sigma_2:\mathbb{H}\times S\longrightarrow S$, $\sigma_2(\psi,\mu)=\mu'$, where $[,]_{\mu'}-[,]_{\mu}$ determines $\psi$. The above discussion shows that the map is well-defined. It is clear that the action is transitive.

Let us consider the relationship between the two actions $\sigma_1$ and $\sigma_2$ on $S$.
\begin{prop}
Let $\lambda$ be a deformation of the Leibniz algebra $L$ with base $A$ and let 
$$0\longrightarrow
\mathbb{K}\stackrel{i}{\longrightarrow} B
\stackrel{p}{\longrightarrow}A\longrightarrow 0$$
be a given extension of $A$.
If $u:B \longrightarrow B$ is an automorphism of this extension 
which corresponds to an element $h \in H_{Harr} ^1(A;\mathbb{K})=TA$, then for any deformation $\mu$ of $L$ with base $B$, such that $p_{*}\mu=\lambda$, the difference $[,]_{u_{*}\mu}-[,]_{\mu}$ is a cocycle in  the cohomology class $d\lambda(h)$. This means that the operation $\sigma_1$ and $\sigma_2$ on $S$ are related to each other by the differential $d\lambda: TA \longrightarrow \mathbb{H}$.
\end{prop}
\begin{proof}
Recall that $\sigma_1:H_{Harr} ^1(A;\mathbb{K})\times S \longrightarrow S,~  \sigma_1(\phi,\mu)=u_{*}\mu$  
where $u$ is the automorphism in $\mathcal{A}$ corresponding to $h =[\phi]\in H_{Harr} ^1(A;\mathbb{K})$, and
$$ \sigma_2:\mathbb{H}\times S \longrightarrow S,~ \sigma_2(\psi,\mu)=\mu^\prime$$
where $[,]_{\mu^\prime}-[,]_\mu$  determines $\psi \in \mathbb{H}$.

We need to show that the $2$-cocycle determined by the difference $[,]_{u_{*}\mu}-[,]_{\mu}$ can be represented by $d \lambda(h)$. One can easily prove this by choosing a basis $\{\tilde{m}_i\}_{1\leq i \leq l}$ of $\mathfrak{M}/\mathfrak{M}^2$ and writing out $[,]_{u_{*}\mu}-[,]_{\mu}$ in terms of $\tilde{m}_i$. Namely, using the notations of Section \ref{UID}, one gets 
$$[,]_{u_{*}\mu}-[,]_{\mu}=\sum_{i=1}^l h(\tilde{m}_i)\otimes \psi_i(l_1,l_2) .$$
Thus the cocycle determined by this difference is $\sum_{i=1}^l h(\tilde{m}_i) \psi_i(l_1,l_2) $. On the other hand for the dual basis $\{\tilde{\xi}_i\}_{1\leq i\leq l}$, if $h=\sum_{i=1}^{l}x_i \tilde{\xi}_i$ and $x_i=h(\tilde{m}_i)$, then $d\lambda(h)=a_{\pi_{*}\lambda}(h)=\sum_{i=1}^{l}h(\tilde{m}_i)[\psi_i]$.
This completes the proof.
\end{proof}

\begin{cor}\label{d is onto}
Suppose that for a deformation $\lambda$ of the Leibniz algebra $L$ with base $A$, the differential $d\lambda: TA \longrightarrow \mathbb{H}$ is onto.
Then the group of automorphisms $\mathcal{A}$ of the extension (\ref{1-dimensional ext}) 
operates transitively on the set of equivalence classes of deformations $\mu$ of $L$ with base $B$ such that $p_{*}\mu=\lambda$. In other words, if $\mu$ exists, it is unique up to an isomorphism and an automorphism of this extension. 
\end{cor}

Suppose now that $M$ is a finite dimensional $A$-module satisfying the condition $\mathfrak{M}M=0$, where $\mathfrak{M}$ is the maximal ideal in $A$. The previous results can  be generalized from the $1$-dimensional extension (\ref{1-dimensional ext}) to a more general extension
$$0\longrightarrow
{M}\stackrel{i}{\longrightarrow} B
\stackrel{p}{\longrightarrow}A\longrightarrow 0.$$

If we try to extend a deformation with base $A$ to a deformation with base $B$, as in the beginning of the this section, then an  analogous computation yields 
$$ \tilde{\phi}:(\frac{B\otimes L}{ker~(E)})^{\otimes 3} \longrightarrow ker~(P)=im (I)\cong M\otimes L.$$
It will give rise to $\bar{\phi} \in CL^3(L;M\otimes L)$ with the cohomology class
 $$[\bar{\phi}]\in HL^3(L;M \otimes L)=M \otimes HL^3(L;L).$$
The obstruction map for this extension is  
$$ \theta_{\lambda}:H_{Harr}^2(A;M) \longrightarrow M \otimes HL^3(L;L) ~\mbox{defined by}~~ \theta_{\lambda}([f])=[\bar{\phi}] .$$
Then, as in the case of $1$-dimensional extension, we have the following.
\begin{prop}\label{for any extension}
 Let $\lambda$ be a deformation of a Leibniz algebra $L$ with base $(A,\mathfrak{M})$ and let $M$ be a finite dimensional $A$-module with $\mathfrak{M}M=0$. Consider an extension $B$ of $A$
 $$0\longrightarrow
{M}\stackrel{i}{\longrightarrow} B
\stackrel{p}{\longrightarrow}A\longrightarrow 0$$
corresponding to some $[f] \in H^{2}_{Harr}(A;M)$. A deformation $\mu$ of $L$ with base $B$ such that $p_{*}\mu=\lambda$ exists if and only if the obstruction $\theta_{\lambda}([f])=0$.
If $d\lambda: TA \longrightarrow \mathbb{H}$ is onto, then the deformation $\mu$, if it exists, is unique up to an isomorphism and an automorphism of the above extension. 
\end{prop}
We end this section with the following naturality property of the obstruction map.
 
\begin{prop}\label{obstructions are same}
Suppose  $A_1$ and $A_2$ are finite dimensional local algebras with augmentations $\varepsilon_1$ and $\varepsilon_2$, respectively.  Let $\phi : A_{2} \longrightarrow A_{1}$ be an algebra  homomorphism with $\phi(1)=1$ and  ${\varepsilon}_1 \circ \phi=\varepsilon_2$. Suppose $\lambda_2$ is a deformation of a Leibniz algebra $L$ with base $A_2$ and $\lambda_1 = \phi_{*} \lambda_2$ is the push-out via $\phi$. Then  the following diagram commutes.
\begin{figure}[htb]
  \begin{center}
    \includegraphics[width=4.5cm]{diag2.epsi}
  \end{center}
  \caption{}
\end{figure}
\end {prop}
\begin{proof}
Let $[f_{A_1}] \in H_{Harr}^2({A_1};\mathbb{K})$ and $[f_{A_2}]=\phi^{*}([f_{A_1}]) \in H_{Harr}^2(A_{2};\mathbb{K})$ correspond to the classes of $1$- dimensional extensions of $A_1$  and $A_2$, represented by 
$$0\longrightarrow
\mathbb{K}\stackrel{i_k}{\longrightarrow} {A}^\prime_{k}
\stackrel{p_k}{\longrightarrow}A_{k}\longrightarrow 0,~~k=1,2.$$
Fix some sections $q_k:A_{k} \longrightarrow A^\prime_{k}$ of $p_k$ for $k=1,2$. Then, as in (\ref{induced by the section}), we get  $\mathbb{K}$-module isomorphisms $A^\prime_{k} \cong A_k \oplus \mathbb{K}$~. Let $(b,x)_{q_k} $ denote the inverse of $(b,x)$  under the above isomorphisms. The algebra structures on $A^\prime_k$ are determined as in (\ref{algebra structure}). Define $\psi:A^\prime_{2} \cong (A_{2} \oplus \mathbb{K}) \longrightarrow A^\prime_{1} \cong (A_{1} \oplus \mathbb{K})$ by $\psi((a,x)_{q_2})= (\phi(a),x)_{q_1}$  for $(a,x)_{q_2} \in A^\prime_{2}$. It is clear that $\psi$ is a $\mathbb{K}$-algebra homomorphism. 
Thus we get a morphism between the two extensions given by
\begin{figure}[htb]
\begin{center}
\includegraphics[width=7.5cm]{obs2.epsi}
\end{center}
\caption{}
\end{figure}

Let $I_{k}=i_{k}\otimes id$, $P_{k}=p_k\otimes id$ and $E_{k}=\hat{\varepsilon}_{k}\otimes id$ , where $\hat{\varepsilon}_k=\varepsilon_{k} \circ p_1$ for $k=1,2$.
Suppose $\mathfrak{M}_{A_k}$ is the unique maximal ideal in $A_k$. Then $\mathfrak{M}_{A^\prime_{k}}=p^{-1}_k(\mathfrak{M}_{A_k})$ is the unique maximal ideal of $A^\prime_{k}$. Denote by $\{m_{ki}\}_{1\leq i \leq r_k}$ a basis of $\mathfrak{M}_{A_k}$ and $\{n_{ki}\}_{1\leq i \leq r_{k}+1}$ a basis of $\mathfrak{M}_{A^\prime_{k}}$ for $k=1,2$. Here $n_{kj}=(m_{kj},0)_{q_k}$ for $1\leq j\leq r_k$  and $n_{k(r_{k}+1)}=(0,1)_{q_k}$.
By (\ref{expression}), the Leibniz bracket on $A_{2}\otimes L$ is defined by
$$ [1\otimes l_1,1\otimes l_2]_{\lambda_2}=1\otimes [l_1,l_2]+\sum_{i=1}^{r_2} m_{2i} \otimes \psi^{2}_i(l_1,l_2)~~\mbox{for} ~ l_1, l_2 \in  L$$ and $\psi^{2}_{i}=\alpha_{\lambda_2 \xi_{2i}}$, where $\{\xi_{2 i}\}$ is the dual basis of $\{m_{2 i}\}$. 
Let  $\phi(m_{2i})=\sum_{j=1}^{r_1}c_{i,j} m_{1j}$ where $c_{i,j} \in \mathbb{K}$ for $1\leq i \leq r_2 ~~\mbox{and}~1\leq j \leq r_1$. Then the push-out $\lambda_1=\phi_{*}\lambda_2$ on $A_{1}\otimes L$ is defined by 
\begin{equation*}
\begin{split}
[1\otimes l_1,1\otimes l_2]_{\lambda_1}
=&~~1\otimes [l_1,l_2]+\sum_{i=1}^{r_2}(\sum_{j=1}^{r_1} c_{i,j} m_{1j})\otimes \psi^{2}_i(l_1,l_2) \\
=&~~1\otimes [l_1,l_2]+\sum_{j=1}^{r_1} m_{1j}\otimes \psi^{1}_j(l_1,l_2)~~\mbox{for}~ l_1, l_2 \in  L.
\end{split}
\end{equation*}
Here $\psi^{1}_j \in CL^2(L;L)$ is defined by
$ \psi^{1}_j(l_1,l_2)=\sum_{i=1}^{r_2} c_{i,j}\psi^{2}_i(l_1,l_2)~~\mbox{for}~ l_1,l_2 \in L.$
For any $2$-cochain $\chi \in CL^2(L;L)$, let
$\{ ,\}_{k}:(A^\prime_{k} \otimes L)^{\otimes 2} \longrightarrow A^\prime_{k} \otimes L$ be the $A^\prime_{k}$-bilinear operation on $A^\prime_{k} \otimes L$ lifting $\lambda_{k}$, defined by
$$\{  1 \otimes l_1,1 \otimes l_2 \}_{k}=1 \otimes [l_1,l_2]+\sum_{j=1}^{r_k} n_{kj} \otimes \psi^{k}_j(l_1,l_2) + n_{k(r_{k}+1)}\chi(l_1,l_2)$$
$\mbox{for}~~ k=1,2 ~~\mbox{and}~l_1, l_2 \in  L$. We know that $\{,\}_{k}$ satisfies properties ({\it i}) and ({\it ii}) of (\ref{condition for extension}).

We claim that $\psi \otimes id$ preserves the liftings. It is enough to show that $(\psi \otimes id)( \{ 1 \otimes l_1,1 \otimes l_2\}_{2})=\{\psi \otimes id~(1 \otimes l_1),\psi \otimes id~(1 \otimes l_2)\}_{1} $ for $l_1,l_2 \in L$. Now
\begin{equation*}
\begin{split}
&(\psi \otimes id)( \{ 1 \otimes l_1,1 \otimes l_2\}_{2})\\
=&~\psi(1) \otimes [l_1,l_2]+\sum_{j=1}^{r_2}  \psi(1)\psi(n_{2j}) \otimes \psi^{2}_j(l_1,l_2) 
+ \psi(1) \psi(n_{2(r_{2}+1)}) \otimes \chi(l_1,l_2) \\
=&~1 \otimes [l_1,l_2]+\sum_{j=1}^{r_2} (\sum_{i=1}^{r_1}c_{j,i} m_{1i}) \otimes \psi^{2}_j(l_1,l_2) +  n_{1(r_{1}+1)}\otimes \chi(l_1,l_2)\\
&( \phi(m_{2j})=\sum_{i=1}^{r_1}c_{j,i} m_{1i} ~\mbox{and}~~\psi(n_{2(r_2+1)})=\psi((0,1)_{q_2})=(\phi(0),1)_{q_1}=n_{1(r_1+1)}) \\
=&~\psi(1) \otimes [l_1,l_2]+\sum_{i=1}^{r_1}  \psi(1) m_{1i} \otimes \psi^{1}_i(l_1,l_2) + \psi(1) n_{1(r_1+1)}\otimes \chi(l_1,l_2)\\ 
=&~ \{\psi(1)\otimes l_1,\psi(1)\otimes l_2\}_{1} \\ 
=&~ \{\psi \otimes id~(1 \otimes l_1),\psi \otimes id~(1 \otimes l_2)\}_{1}, ~\mbox{ which proves our claim}.
\end{split}
\end{equation*}
Let $\phi_k$ be defined by $\{,\}_{k}$ as in (\ref{phi}) and $\bar{\phi}_{k}$ the corresponding cocycle as in (\ref{relation of phi}). As $\psi(n_{2(r_2+1)})=n_{1(r_1+1)}$, it follows from the definition of $\phi_{k}$ and the previous claim that $[\bar {\phi_{2}}] =[\bar {\phi_{1}}]$. 
Therefore,
$\theta_{\lambda_{1}}([f_{A_1}])=[\bar{\phi_{1}}] 
=[\bar {\phi_{2}}] 
= \theta _{\lambda_{2}}([f_{A_2}]) 
= \theta _{\lambda_{2}} \circ \phi^{*}([f_{A_{1}}]).$
Hence $\theta_{\lambda_{1}}=~\theta _{\lambda_2} \circ \phi^{*}$.
\end{proof}
 
\section {Construction of a Versal Deformation }
In this section we give an explicit construction of versal deformation of a given Leibniz algebra following \cite{FiFu}.

Consider the Leibniz algebra $L$ with $dim(\mathbb{H})<\infty$.
Set $C_0=\mathbb{K}$ and 
$C_1=\mathbb{K}\oplus \mathbb{H}^\prime$.
Consider the extension
$$0\longrightarrow
 \mathbb{H}^\prime \stackrel{i}{\longrightarrow} C_1
\stackrel{p}{\longrightarrow} C_0 \longrightarrow 0,$$ 
where the multiplication in $C_1$ is defined by 
$$(k_1,h_1) \cdot (k_2,h_2) =(k_1 k_2~ ,~ k_1 h_2+k_2 h_1)~~\mbox{for}~ (k_1,h_1) , (k_2,h_2)\in C_1.$$
Let $\eta_1$ be the universal infinitesimal deformation with base $C_1$ as constructed in Section \ref{UID}.
We  proceed by induction.  Suppose for some $k\geq 1$ we have constructed a finite dimensional local algebra $C_k$ and a  deformation $\eta_k$ of $L$ with base $C_k$.   
Let 
$$\mu:H_{Harr}^2(C_k;\mathbb{K})\longrightarrow (Ch_{2} (C_k))^\prime$$
be a homomorphism sending a cohomology class to a cocycle representing the class. Let
$$f_{C_k}:Ch_{2} (C_k) \longrightarrow H_{Harr}^2(C_k;\mathbb{K})^\prime$$
 be the dual of $\mu$. By Proposition \ref{cohomology class corresponds to extension} ({\it ii}) we have the following extension of $C_k$:  
\begin{equation}\label{universal extension}
0\longrightarrow
 H_{Harr}^2(C_k;\mathbb{K})^\prime \stackrel{\bar {i}_{k+1}}{\longrightarrow} {\bar C}_{k+1}
\stackrel{\bar{p}_{k+1}}{\longrightarrow} C_k \longrightarrow 0.
\end{equation}
The corresponding obstruction $\theta([f_{C_k}]) \in H_{Harr}^2(C_k;\mathbb{K})^\prime \otimes HL^3(L;L)$ gives a linear map   
$\omega_k:H_{Harr}^2(C_k;\mathbb{K}) \longrightarrow HL^3(L;L)$
with the dual map  
$${\omega_k}^\prime:HL^3(L;L)^\prime \longrightarrow H_{Harr}^2(C_k;\mathbb{K})^\prime .$$
We have an induced extension 
$$ 0\longrightarrow coker (\omega'_{k})\longrightarrow \bar{C}_{k+1}/\bar{i}_{k+1}\circ \omega'_{k}(HL^3(L;L)')\longrightarrow C_k \longrightarrow 0.$$
Since $coker (\omega'_k)\cong (ker (\omega_k))^\prime$,
it yields an extension 
\begin{equation}\label{yields an extension}
0\longrightarrow (ker(\omega_k))^\prime \stackrel{i_{k+1}}{\longrightarrow} C_{k+1}
\stackrel{p_{k+1}}{\longrightarrow} C_k \longrightarrow 0
\end{equation}
where $C_{k+1}= { \bar{C}_{k+1}}/{\bar{i}_{k+1}\circ~ \omega_{k}^\prime (HL^3(L;L)')}$ and  $i_{k+1}$, $p_{k+1}$ are the mappings induced by $\bar{i}_{k+1}$ and $\bar{p}_{k+1}$, respectively.
Observe that the algebra $C_k$ is also local.  Since $C_k$ is finite dimensional, the cohomology group $H_{Harr}^2(C_k;\mathbb{K})$ is also finite dimensional and hence $C_{k+1}$ is finite dimensional as well. 
\begin{rem}\label{specific extension}
It follows from Proposition \ref{coefficients in M} that the specific extension (\ref{universal extension}) has the following ``universality property''. For any $C_k$-module $M$ with $\mathfrak{M}M=0$, (\ref{universal extension}) admits a unique morphism into an arbitrary extension of $C_k$:
$$0\longrightarrow M \longrightarrow B \longrightarrow C_k \longrightarrow 0.$$ 
\end{rem}

\begin{prop}
The deformation $\eta_k$ with base $C_{k}$ of a Leibniz algebra $L$ admits an extension to a deformation with base $C_{k+1}$, which is unique up to an isomorphism and an automorphism of the extension
$$0\longrightarrow
 (ker(\omega_k))^\prime \stackrel{i_{k+1}}{\longrightarrow} C_{k+1}
\stackrel{p_{k+1}}{\longrightarrow} C_k \longrightarrow 0.$$
\end{prop}
\begin{proof}
From the above construction of the extension (\ref{yields an extension})
it is clear that the corresponding obstruction map is the restriction of $\omega_{k}$,
$$\omega_{k}|_{ker(\omega_k)}:ker(\omega_k)\longrightarrow HL^3(L;L).$$
Hence, it is the zero map.  Thus the result follows from Proposition \ref{for any extension}.
\end{proof}
By induction, the above process yields a sequence of finite dimensional local algebras $C_{k}$ and deformations $\eta_{k}$ of the Leibniz algebra $L$ with base $C_{k}$ 
$$ \mathbb{K} \stackrel{p_{1}}{\longleftarrow} C_{1} \stackrel{p_{2}}{\longleftarrow} C_{2}\stackrel{p_{3}}{\longleftarrow} \ldots \ldots \stackrel{p_{k}}{\longleftarrow} C_{k}\stackrel{p_{k+1}}{\longleftarrow} C_{k+1}\ldots$$
such that $ {p_{k+1}}_{*} \eta_{k+1}=\eta_{k}$.
Thus by taking the projective limit we obtain a formal deformation $\eta$ of $L$ with base $C=\ilim_{k\rightarrow
\infty} C_{k}$.

Next, we give an algebraic description of the base $C$ of the versal deformation. For that we need the following Proposition  from \cite{H}.
\begin{prop}\label{AMI}
Let $A=\mathbb{K}[x_1,x_2,\ldots,x_n]$ be the polynomial algebra, and let $\mathfrak{M}$ be the ideal of polynomials without constant terms.\\
$(a)$ If an ideal $I$ of $A$ is contained  in $\mathfrak{M}^2$, then $H_{Harr}^2(A/I;\mathbb{K}) \cong (I/{\mathfrak{M}I})'$. \\
$(b)$ There is an extension for $B=A/I$: 
$$0\longrightarrow {I/{\mathfrak{M}I}}    \stackrel{i}{\longrightarrow} A/{\mathfrak{M}I}
\stackrel{p}{\longrightarrow} A/I \longrightarrow 0$$
where $i$ and $p$ are induced by the inclusions $I \hookrightarrow A$ and $ \mathfrak{M}I \hookrightarrow I$.
\end{prop}

Suppose $dim (\mathbb{H})=n$. Let $\{h_i\}_{1\leq i \leq n}$ be a basis of $\mathbb{H}$ and $\{g_i\}_{1\leq i \leq n}$ be the corresponding dual basis. Let $\mathbb{K}[[\mathbb{H}']]$ denote the formal power series ring $\mathbb{K}[[g_1,\ldots,g_n ]]$ in $n$ variables $g_1,\ldots,g_n$ over $\mathbb{K}$~. Now a typical element in $\mathbb{K}[[\mathbb{H}^\prime]]$ is of the form $$\sum_{i=0}^{\infty}a_i f_i(g_1,\ldots,g_n)=a_0+a_1f_1(g_1,\ldots,g_n)+a_2f_2(g_1,\ldots,g_n)+\ldots $$
where  $a_i \in \mathbb{K}$ and $f_{i}$ is a monomial of degree $i$ in $n$-variables $g_1,\ldots,g_n$ for $i=0,1,2,\ldots$ . Let $\mathfrak{M}$ denote the unique maximal ideal in $\mathbb{K}[[\mathbb{H}^\prime]]$, consisting of all elements in $\mathbb{K}[[\mathbb{H}']]$ with constant term being equal to zero. 
\begin{prop}
For the local algebra $C_k$ we have $C_k \cong \mathbb{K}[[\mathbb{H}^\prime]]/I_{k}$ for some ideal $I_k$, satisfying  $\mathfrak{M}^2=I_1\supset I_2\supset \ldots \supset I_{k}\supset \mathfrak{M}^{k+1}$.
\end{prop}
\begin{proof}
 By construction,
$C_1=\mathbb{K}\oplus \mathbb{H}^\prime \cong \mathbb{K}[[\mathbb{H}^\prime]]/\mathfrak{M}^2.$
Suppose we already know that 
$C_{k} \cong \mathbb{K}[[\mathbb{H}^\prime]]/I_{k}$
where $\mathfrak{M}^2 \supset I_{k}\supset \mathfrak{M}^{k+1}$. Then  by specifying $A=\mathbb{K}[[\mathbb{H}']]$ and $I=I_k$ in Proposition \ref{AMI}, we get  $\bar{C}_{k+1} \cong \mathbb{K}[[\mathbb{H}^\prime]]/\mathfrak{M} I_{k}.$
In the previous construction, $C_{k+1}$ is the quotient of $\bar{C}_{k+1}$ by an ideal contained in $I_{k}/\mathfrak{M} I_{k}\subset \mathfrak{M}^2/\mathfrak{M} I_{k}$. 
Hence $C_{k+1} \cong \mathbb{K}[[\mathbb{H}^\prime]]/I_{k+1}$
where $\mathfrak{M}^2\supset I_{k+1}\supset \mathfrak{M}I_{k}\supset \mathfrak{M}^{k+2}$. The proof is now complete by induction.
\end{proof}

\begin{cor}\label{about differential}
For $k\geq 2$ the projection $p_{k}:C_{k}\longrightarrow C_{k-1}$ induces an isomorphism $TC_{k} \longrightarrow TC_{k-1}$. In particular, for every $k\geq 1$, $TC_{k}\cong TC_{1}=\mathbb{H}$.  Moreover, under the above identification of $TC_k$ with $\mathbb{H}$, the differential $d\eta_{k}:TC_{k}\longrightarrow \mathbb{H}$ is the identity map. 
\end{cor}
\begin{proof}
 We have $C_{0}= \mathbb{K}$ ; $C_{1}=\mathbb{K}\oplus \mathbb{H}^\prime \cong \mathbb{K}[[\mathbb{H}^\prime]]/\mathfrak{M}^2$ and for $k \geq 2$, 
 $C_k=\mathbb{K}[[\mathbb{H}^\prime]]/I_{k}$ 
 where $\mathfrak{M}^2=I_1\supset I_2\supset \ldots \supset I_{k}\supset \mathfrak{M}^{k+1}$.
  The projection $p_{k}:C_{k}\longrightarrow C_{k-1}$ is given by $p_{k}(f+I_k)=f+ I_{k-1}~~ \mbox{for}~f \in C_k ~\mbox{and}~k\geq 1$. 
The map $p_{k}$  gives rise to a surjective linear map $ \mathfrak{M}/I_k \longrightarrow \mathfrak{M}/I_{k -1}$. Taking the quotient map $\mathfrak{M}/I_{k -1}\longrightarrow \frac{\mathfrak{M}/I_{k -1}}{\mathfrak{M}^2/I_{k -1}} $,  we get an epimorphism $\mathfrak{M}/I_k \longrightarrow \frac{\mathfrak{M}/I_{k -1}}{\mathfrak{M}^2/I_{k -1}} $ with kernel $\mathfrak{M}^2/I_{k}$ which corresponds to an isomorphism
$$ \frac{\mathfrak{M}/I_{k }}{\mathfrak{M}^2/I_{k }}\longrightarrow \frac{\mathfrak{M}/I_{k -1}}{\mathfrak{M}^2/I_{k -1}}.$$
As a result we get an isomorphism
$$( \frac{\mathfrak{M}/I_{k }}{\mathfrak{M}^2/I_{k }} )^\prime =TC_{k}\longrightarrow TC_{k-1}=(\frac{\mathfrak{M}/I_{k -1}}{\mathfrak{M}^2/I_{k -1}})^\prime. $$
Observe that for any $k\geq 1$, $TC_k=( \frac{\mathfrak{M}/I_{k }}{\mathfrak{M}^2/I_{k }} )^\prime =(\frac{\mathfrak{M}}{\mathfrak{M}2})^\prime \cong TC_1$. 
On the other hand, since $C_1=\mathbb{K}\oplus \mathbb{H}^\prime$ with maximal ideal $\mathbb{H}'$ and $(\mathbb{H}')^2=0$. Hence $TC_1=(\mathbb{H}')'=\mathbb{H}$. The last assertion follows from the definition of the differential.
\end{proof}
\begin{prop}
The complete local algebra $C=\ilim_{k\rightarrow \infty} C_{k}$ can be described as $C \cong \mathbb{K}[[\mathbb{H}^\prime]]/I$, where $I$ is an ideal contained in $\mathfrak{M}^2$.
\end{prop}
\begin{proof}
Consider the map 
$$\phi:\mathbb{K}[[\mathbb{H}^\prime]]\longrightarrow C_k=\mathbb{K}[[\mathbb{H}^\prime]]/I_k~~\mbox{defined by}~\phi(f)=f+I_k~~\mbox{for}~~f\in \mathbb{K}[[\mathbb{H}^\prime]].$$
Since $I_k \supset \mathfrak{M}^{k+1}$, the map $\phi$ induces an epimorphism 
$$\phi_{k}:\mathbb{K}[[\mathbb{H}^\prime]]/\mathfrak{M}^{k+1}\longrightarrow C_k~~\mbox{for each}~k\geq 1.$$
In the limit we get an epimorphism
$$ \mathbb{K}[[\mathbb{H}']]=\ilim_{k\rightarrow
\infty} \mathbb{K}[[\mathbb{H}^\prime]]/\mathfrak{M}^{k+1} \longrightarrow \ilim_{k\rightarrow
\infty} C_{k}.$$
Therefore $C\cong \mathbb{K}[[\mathbb{H}^\prime]]/I$ where $I=\bigcap_k I_k $ is the kernel of the epimorphism.
\end{proof}

Finally we prove the versality property of the constructed deformation $\eta$ with base $C$. For this we use the following standard lemma.
\begin{lem}\label{splits into}
Suppose $0\longrightarrow
{M_s}\stackrel{i}{\longrightarrow} B_s
\stackrel{p}{\longrightarrow}A\longrightarrow 0  $ is an $s$-dimensional extension of $A$. Then there exists an $(s-1)$-dimensional extension $$ 0\longrightarrow
{M_{s-1}}\stackrel{\bar{i}}{\longrightarrow} B_{s-1}
\stackrel{\bar{p}}{\longrightarrow}A\longrightarrow 0$$ of $A$ and a $1$-dimensional extension $$0\longrightarrow
{\mathbb{K}}\stackrel{i'}{\longrightarrow} B_s
\stackrel{p'}{\longrightarrow}B_{s-1}\longrightarrow 0.$$
\end{lem}
\begin{thm}
Let $L$ be a Leibniz algebra with $dim(\mathbb{H})<\infty$. Then the formal deformation $\eta$ with base $C$ constructed above is a versal deformation of $L$.
\end{thm}
\begin{proof}
Suppose $dim(\mathbb{H})=n$. Let $\{h_i \}_{1\leq i \leq n}$ be a basis of  
$\mathbb{H}$ and $\{g_i \}_{1\leq i \leq n}$ the corresponding dual basis of $\mathbb{H}^\prime$. Let $A$ be a complete local algebra with maximal ideal $\mathfrak{M}$ and let $\lambda$ be a formal deformation of $L$ with base $A$. We want to find a $\mathbb{K}$-algebra homomorphism $\phi:C \longrightarrow A$ such that $\phi _{*} \eta=\lambda$. Denote $A_0=A/\mathfrak{M}\cong \mathbb{K}~;~A_1=A/\mathfrak{M}^2\cong \mathbb{K}\oplus (TA)^\prime$. Since $A$ is complete, we have $A=\ilim_{k\rightarrow
\infty} A/{\mathfrak{M}}^k $. Moreover, for each $k$ we have the following  finite dimensional extension
$$0\longrightarrow \frac {{\mathfrak{M}}^k}{{\mathfrak{M}}^{k+1}}\longrightarrow \frac{A}{{\mathfrak{M}}^{k+1}}\longrightarrow \frac{A}{{\mathfrak{M}}^{k}}\longrightarrow 0,$$
because $dim (\frac {{\mathfrak{M}}^k}{{\mathfrak{M}}^{k+1}})<\infty$.

Let $dim (\frac {{\mathfrak{M}}^k}{{\mathfrak{M}}^{k+1}})=n_{k-1}$. A repeated application of Lemma \ref{splits into} to the extension $$0\longrightarrow \frac {{\mathfrak{M}}^2}{{\mathfrak{M}}^3}\longrightarrow \frac{A}{\mathfrak{M}^3}\longrightarrow \frac{A}{\mathfrak{M}^2}=A_1\longrightarrow 0$$
yields $n_1$ number of $1$-dimensional extensions as follows.
\begin{equation*}
\begin{split}
& 0\longrightarrow \mathbb{K}\longrightarrow A_2 \longrightarrow A_1 \longrightarrow 0\\
&0\longrightarrow \mathbb{K}\longrightarrow A_3 \longrightarrow A_2 \longrightarrow 0\\
&~~~~~~~~~~~~~~~~~~~\vdots \\
&0\longrightarrow \mathbb{K}\longrightarrow A_{n_1 +1}=\frac{A}{\mathfrak{M}^3} \longrightarrow A_{n_1} \longrightarrow 0.
\end{split}
\end{equation*}
Similarly, the extension 
$$0\longrightarrow \frac {{\mathfrak{M}}^3}{{\mathfrak{M}}^4}\longrightarrow \frac{A}{\mathfrak{M}^4}\longrightarrow \frac{A}{\mathfrak{M}^3}=A_{n_1+1}\longrightarrow 0$$
splits into $n_2$ number of $1$-dimensional extensions and so on.
Thus we get a sequence of $1$- dimensional extensions
$$0\longrightarrow
\mathbb{K}\stackrel{j_{k+1}}{\longrightarrow} A_{k+1}
\stackrel{q_{k+1}}{\longrightarrow}A_{k}\longrightarrow 0~~~;~k\geq 1.$$
Since $A=\ilim_{k\rightarrow\infty} A/{\mathfrak{M}}^k $, it is clear that $A=\ilim_{k\rightarrow
\infty} A_{k} $.
 Let $Q_{k}:A \longrightarrow A_{k}$ be the projection map for the inverse system $\{A_k,q_k\}_{k\geq 1}$ with the limit $A$,  where $Q_1:A\longrightarrow A_1=A/\mathfrak{M}^2$ is the natural projection.  Let ${Q_{k}}_{*}\lambda=\lambda_{k}$, then $\lambda_{k}$ is a deformation of $L$ with base $A_{k}$. Thus $\lambda_k= {Q_{k}}_*\lambda =(q_{k+1}\circ Q_{k+1})_{*}\lambda ={q_{k+1}}_{*}\lambda_{k+1}$.  Now we will construct inductively homomorphisms $\phi_j:C_j\longrightarrow A_j$ for $j=1,2\ldots$, compatible with the corresponding projections $C_{j+1}\longrightarrow C_j$ and $A_{j+1}\longrightarrow A_{j}$, along with the conditions ${\phi_j}_{*} \eta_j \cong\lambda_j$. Define $$\phi_1:C_1 \longrightarrow A_1~~\mbox{as}~~id \oplus (d\lambda)^\prime :\mathbb{K}\oplus \mathbb{H}^\prime \longrightarrow \mathbb{K}\oplus (TA)^\prime.$$
From Proposition \ref{couniversal} we have  ${\phi_1}_{*}\eta_1 \cong \lambda_1$.  

Suppose we have constructed a $\mathbb{K}$-algebra homomorphism $\phi_k: C_k\longrightarrow A_k$    with ${\phi_k}_*\eta_k \cong \lambda_k$.  Consider the homomorphism  $\phi_{k}^{*}:H_{Harr}^2(A_k;\mathbb{K})\longrightarrow H_{Harr}^2(C_k;\mathbb{K})$ induced by $\phi_{k}$. Let $$ 0 \longrightarrow \mathbb{K}\stackrel{i_{k+1}}{\longrightarrow} B \stackrel{p_{k+1}}{\longrightarrow} C_k \longrightarrow0$$ represent the image under $\phi^{*}_{k}$ of the isomorphism class of extension 
$$0 \longrightarrow \mathbb{K}\stackrel{j_{k+1}}{\longrightarrow} A_{k+1}\stackrel{q_{k+1}}{\longrightarrow} A_{k} \longrightarrow 0$$ (see Proposition \ref{cohomology class corresponds to extension}). Then we have the following commutative diagram
\begin{figure}[htb]
  \begin{center}
    \includegraphics[width=7.5cm]{tab1.epsi}
  \end{center}
 \caption{} 
\end{figure}

where $\psi$ is given by $\psi((x,k)_q)=(\phi_k(x),k)_{q'}$ for some fixed sections $q$ and $q'$ of $p_{k+1}$ and $q_{k+1}$ respectively. Observe that by Proposition \ref{obstructions are same} the obstructions in extending $\lambda_{k}$ to the base $A_{k+1}$ and that of $\eta_k$ to the base $B$ coincide. Since $\lambda_k$ has an extension $\lambda_{k+1}$, the corresponding obstruction is zero. Hence there exists a deformation $\xi$ of $L$ with base $B$ which extends $\eta_{k}$ with base $C_{k}$  such that  ${\psi}_{*} \xi=\lambda_{k+1}$. By Remark \ref{specific extension} we get the following unique morphism of extensions.
\begin{figure}[htb]
  \begin{center}
    \includegraphics[width=7.5cm]{table2.epsi}
  \end{center}
 \caption{} 
\end{figure}

Since the deformation $\eta_{k}$ has been extended to $B$, the obstruction map
$$\omega_k: H_{Harr}^2(C_{k};\mathbb{K})\longrightarrow HL^3(L;L)$$ is zero. Therefore the composition $\tau^\prime \circ \omega^\prime _{k}: HL^3(L;L)'\longrightarrow \mathbb{K}$ is zero. So  $\tau^\prime$ will induce a linear map 
$\tau:H_{Harr}^2(C_{k};\mathbb{K})'/\omega^\prime _{k}(HL^3(L;L)')\longrightarrow \mathbb{K}$.  
Also the map $\bar{\chi}:\bar{C}_{k+1}\longrightarrow B$ will induce a linear map $\chi:C_{k+1}=\bar{C}_{k+1}/\bar{i}_{k+1}\circ \omega^\prime _{k}(HL^3(L;L)')\longrightarrow B.$
Since $coker(\omega^\prime _{k})\cong (ker (\omega_{k}))^\prime$, the last diagram yields the following commutative diagram.
\begin{figure}[htb]
  \begin{center}
    \includegraphics[width=7.5cm]{tab4.epsi}
  \end{center}
 \caption{}
\end{figure}

By Corollary \ref{about differential}, the differential 
$$d\eta_k:TC_k \longrightarrow \mathbb{H}$$
is onto, so by Corollary \ref{d is onto}, the deformations $\chi_*\eta_{k+1}$ and $\xi$ are related by some automorphism $u:B\longrightarrow B $ of the extension 
$$ 0\longrightarrow \mathbb{K}\longrightarrow B \longrightarrow C_k \longrightarrow0$$ with $u_{*}(\chi_{*}\eta_{k+1})=\xi$.
Now set  $ \phi_{k+1}=(\psi \circ u \circ \chi):C_{k+1}\longrightarrow A_{k+1}$, where $\psi$ is as in Figure 3 .
Then 
\begin{equation*}
\begin{split}
{\phi_{k+1}}_{*}\eta_{k+1}=&~ \psi_{*}\circ u_{*}\circ \chi_{*} \eta_{k+1}
=~\psi_{*} \xi 
=~\lambda_{k+1}.
\end{split}
\end{equation*}
Thus by induction we get a sequence of homomorphisms $\phi_k:C_k\longrightarrow A_k$ with ${\phi_{k}}_*\eta_k=\lambda_k$. Consequently, taking the limit, we find a homomorphism  $\phi:C\longrightarrow A$ such that $\phi_* \eta=\lambda$.
 If $\mathfrak{M}^2=0$, then the uniqueness of $\phi$ follows from the corresponding property in Proposition \ref{couniversal}.
\end{proof}

\section{Conclusions}
In this work we gave a constructive method for Leibniz algebras for the solution of the main deformation question, suitable for specific computations. The main feature of this method is that it completely describes all non-equivalent deformations $-$ a problem which did not have a satisfactory solution for a long time. For this we had to consider deformations with complete local algebra base, which was necessary for the existence of a versal deformation. The construction presented here is an inductive procedure, which consists of extending the base of deformation at each step. The specific description of the base of the versal deformation is useful for computations.

{\bf Alice Fialowski}\\
E$\ddot{o}$tv$\ddot{o}$s Lor$\acute{a}$nd University, Budapest, Hungary.\\
e-mail: fialowsk@cs.elte.hu

{\bf Ashis Mandal}\\
Indian Statistical Institute, Kolkata, India.\\
e-mail: ashis\_r@isical.ac.in

{\bf Goutam Mukherjee}\\
Indian Statistical Institute, Kolkata, India.\\
e-mail: goutam@isical.ac.in
\end{document}